\documentclass[11pt]{article}
\pdfoutput=1

\usepackage[truedimen,margin=25truemm]{geometry}
\usepackage[numbers,compress]{natbib}

\usepackage[utf8]{inputenc} 
\usepackage[T1]{fontenc}    
\usepackage{hyperref}       
\usepackage{url}            
\usepackage{booktabs}       
\usepackage{amsfonts}       

\usepackage{natbib}
\usepackage{algorithm,algorithmic}
\usepackage{amsthm}
\usepackage{amsmath}
\usepackage{amssymb,bbm}
\usepackage{bm}
\usepackage{dsfont}
\usepackage{graphicx}
\usepackage{nicefrac}       
\usepackage{microtype}      
\usepackage{multirow}
\usepackage{makecell}
\usepackage{caption} 

\newtheorem{thm}{Theorem}
\newtheorem{lem}[thm]{Lemma}

\newtheorem{definition}[thm]{Definition}

\newcommand{\diag}{\operatorname{Diag}}

\newcommand{\Gap}{\operatorname{Gap}}
\newcommand{\R}{\mathbb{R}}
\newcommand{\one}{\mathds{1}}
\newcommand{\mat}[1]{\mathbf{#1}}
\renewcommand{\vec}[1]{\bm{#1}}

\usepackage{xcolor}

\newcommand{\changeXSS}[1]{\textcolor{black}{#1}}

\definecolor{ao(english)}{rgb}{0.0, 0.5, 0.0}

\usepackage{tabularx}

\title{Safe Screening for Unbalanced Optimal Transport}

\author{
Xun SU \thanks{Department of Communications and Computer Engineering, School of Fundamental Science and Engineering, WASEDA University, 3-4-1 Okubo, Shinjuku-ku, Tokyo 169-8555, Japan (e-mail: suxun\_opt@asagi.waseda.jp) }
 \and Zhongxi Fang \thanks{Department of Communications and Computer Engineering, School of Fundamental Science and Engineering, WASEDA University, 3-4-1 Okubo, Shinjuku-ku, Tokyo 169-8555, Japan (e-mail: fzx@akane.waseda.jp) }
 \and Hiroyuki Kasai \thanks{Department of Communications and Computer Engineering, School of Fundamental Science and Engineering, WASEDA University, 3-4-1 Okubo, Shinjuku-ku, Tokyo 169-8555, Japan (e-mail: hiroyuki.kasai@waseda.jp)}
 }

\begin{document}

\maketitle

\begin{abstract}
This paper introduces a framework that utilizes the Safe Screening technique to accelerate the optimization process of the Unbalanced Optimal Transport (UOT) problem by proactively identifying and eliminating zero elements in the sparse solutions. We demonstrate the feasibility of applying Safe Screening to the UOT problem with $\ell_2$-penalty and KL-penalty by conducting an analysis of the solution's bounds and considering the local strong convexity of the dual problem. Considering the specific structural characteristics of the UOT in comparison to general Lasso problems on the index matrix, we specifically propose a novel approximate projection, an elliptical safe region construction, and a two-hyperplane relaxation method. These enhancements significantly improve the screening efficiency for the UOT's without altering the algorithm's complexity.
\end{abstract}

\section{INTRODUCTION}
\label{sec:int}

Optimal transport (OT), as a metric, has gained significant attention in the field of machine learning in recent years due to its remarkable ability to capture geometric relationships between data distributions. It has demonstrated impressive achievements in many fields \cite{Courty_PAMI_2017, arjovsky2017wasserstein, Chen_ICLR_2019,Maretic_NIPS_2019}. To overcome the limitation of OT in handling data with unequal quantities, researchers introduced {\it unbalanced optimal transport} (UOT) \cite{refId0} by relaxing the constraints using penalty functions. UOT has been found extensive applications in computational biology \cite{Schiebinger_CELL_2019}, machine learning \cite{Janati_AISTATS_2019}, and deep learning domains \cite{Yang_ICLR_2019}.

However, compared to traditional metrics, the computational burden associated with OT, including UOT, has impeded their widespread adoption on large-scale problems. The state-of-the-art linear programming algorithms suffer from cubic computational complexity and are challenging to parallelize on GPUs \cite{DBLP:journals/orl/Volgenant05}. To address these challenges, Sinkhorn's algorithm \cite{Sinkhorn_1974} has been introduced for entropy-regularized OT and UOT \cite{Cuturi_NIPS_2013,DBLP:journals/moc/ChizatPSV18}, reducing the computational complexity \cite{DBLP:conf/nips/AltschulerWR17, DBLP:conf/icml/PhamLHPB20} . Furthermore, due to the instability of optimization and loss of sparsity in the solutions of Sinkhorn on entropy regularization \cite{DBLP:journals/corr/Schmitzer16, Blondel_AISTATS_2018}, researchers have shifted their attention towards exploring alternative regularization terms. This has led to the proposal of a series of new optimization algorithms \cite{DBLP:conf/nips/GenevayCPB16, DBLP:conf/icml/GuminovDTG21, Nguyen_arXiv_2022}.

Apart from algorithmic advancements, low-rank approximations in Sinkhorn's iterations have accelerated OT and UOT \cite{DBLP:conf/icml/ScetbonCP21, DBLP:conf/nips/AltschulerBRN19}, by further sacrificing solution accuracy to obtain dense approximate solutions. However, considering the inherent sparsity of OT's solutions, transportation mess tend to concentrate on a small number of transportation pairs. Researchers have attempted to reduce the computational burden by early pruning of infeasible pairs or by neglecting small solution elements in advance \cite{doi:10.1137/16M1106018, DBLP:conf/icml/KlicperaLG21} . However, these methods have remained empirical in nature, lacking rigorous theoretical foundations. 

Recently, researchers has revealed the connection between Lasso-like problems \cite{Tibshirani_JRSS_1996,Efron_AM_2004} and UOT \cite{Chapel_NeurIPS_2021}. This finding has motivated us to consider applying the well-known {\it Safe Screening} technique \cite{ghaoui2010safe} in the lasso domain to the UOT. Safe screening can theoretically identify and freeze the zero elements in the solution before computing the sparse optimum, leveraging the sparsity to shrink the dimension of the problem. Moreover, the sparsity of OT and UOT increases along the dimension, which would benefit a lot for large-scale problems. However, the lack of uniform weights caused by the cost matrix in UOT leads to poor performance when using traditional Safe Screening algorithms based on $\ell_1$-norm regularization, even resulting in theoretical degeneration and failure. Moreover, for the commonly used KL-penalized UOT, the existing method requires the index matrix to be invertible, which is not satisfied by the KL-penalized UOT \cite{9414183}.

This paper presents a theoretical framework for applying Safe Screening to UOT. Additionally, we propose a novel projection method and Safe Screening techniques based on the unique sparse structure of the index matrix in UOT. For an overview of our contributions, please refer to \textbf{Table~\ref{Tab:1}}. Our specific contributions can be summarized as follows:

\begin{itemize}
\item This paper presents the first feasible theoretical framework for Safe Screening on UOT, including the dual problem, the projection method, and dual feasible region construction.
\item Leveraging the specific structural characteristics of the index matrix in UOT, we introduce the {\it Shifting Projection} in {\bf Section~\ref{Sec:SP} }as a new approximate projection method that significantly reduces projection errors without increasing computational complexity, thereby improving screening efficiency.
\item For the KL-penalized UOT, we establish the local blockwise strong convexity property through a bound analysis for a new dual feasible region. Based on the analysis, we propose to construct smaller {\it Ellipse-based} safe regions to exploit the anisotropy of the dual variables in the dual space in {\bf Section~\ref{Sec:Ell} }. To the best of our knowledge, this is the first utilization of ellipses to construct safe regions in the Safe Screening community.
\item Considering the sparse correlation of the index matrix, we propose the {\it Cruciform Two-hyperPlane} (CTP) method to further shrink the safe region in {\bf Section~\ref{Sec:CTP} }. Using the Lagrangian method, we obtain closed-form solutions without significantly increasing the computational burden.
\end{itemize}

\begin{table}[t]
\caption{Feasibility of applying Safe Screening to UOT under different penalty functions is examined in the table. The 'Projection' column indicates the selection of the projection method, with 'R' representing Residuals Rescaling and 'S' representing Shifting Projection. 'p' denotes the presence of infeasible instances in the method.}
\label{sample-table}
\begin{tabular}{l||c|c|ccr}
\toprule
Base &\makecell{Abbre-\\viation} &\makecell{Projection \\ {\bf (Section~\ref{Sec:SP}) }} & \makecell{KL-\\
penalty} &\makecell{$\ell_2$-\\penalty} \\
\midrule
Gap-safe \cite{JMLR:v18:16-577} \cite{DBLP:conf/icml/AtamturkG20} & Gap & R& $\times$ & p & \\
& & S & $\times$ & $\surd$ & \\
\midrule
Dynamic Sasvi \cite{Yamada_NIPS_2021}& Sa& R& $\times$  & p &\\
&  &S & $\times$  & $\surd$ &\\
\midrule
Gap-Ellipse {\bf (Section~\ref{Sec:Ell}) }&Ell & S & $\surd$ & $\times$  &\\
\midrule
Cruciform Two-hyperPlane {\bf (Section~\ref{Sec:CTP}) }&Gap-CTP &S & $\surd$ & $\surd$ &\\
        &Sa-CTP &S & $\times$  & $\surd$ &\\
&Ell-CPT &S& $\surd$ & $\times$  &\\
\bottomrule
\end{tabular}
\label{Tab:1}
\end{table}

\section{PRELIMINARIES}
\label{sec:pre}

$\mathbb{R}^n$ denotes $n$-dimensional Euclidean space, and $\mathbb{R}^n_+$ denotes the set of vectors in which all elements are non-negative. $\mathbb{R}^{m \times n}$ represents the set of $n \times m$ matrices. In addition, $\mathbb{R}^{n \times m}_+$ stands for the set of $n \times m$ matrices in which all elements are non-negative. We present vectors as bold lower-case letters $\vec{a},\vec{b},\vec{c},\dots$ and matrices as bold-face upper-case letters $\mat{A},\mat{B},\mat{C},\dots$. The $i$-th element of $\vec{a}$ and the element at the $(i,j)$ position of $\mat{A}$ are stated respectively as $a_i$ and ${A}_{i,j}$. The $i$-th column of $\mat{A}$ is represented as $\vec{a}_i$. In addition, $\one_n \in \mathbb{R}^n$ is the $n$-dimensional vector in which all elements are one. For two matrices of the same size $\mat{A}$ and $\mat{B}$, $\langle \mat{A},\mat{B}\rangle={\rm tr}(\mat{A}^T\mat{B})$ is the Frobenius dot-product. \changeXSS{We use $\|\vec{a}\|_2$, $\|\vec{a}\|_1$, and $\|\vec a\|_{\infty}$ to represent the $\ell_2$-norm, $\ell_1$-norm, and $\ell_{\infty}$ norm of $\vec{a}$}, respectively. $D_\phi$ is the Bregman divergence with the strictly convex and differentiable function $\phi$, i.e., $D_\phi(\vec{a},\vec{b})=\sum_{i} d_\phi(a_i, b_i)=\sum_i [\phi(a_i) - \phi(b_i) - \phi'(b_i)(a_i -b_i)]$. Additionally, we suggest vectorization for $\mat{A} \in \mathbb{R}^{n \times m}$ as lowercase letters $\vec{a} \in \mathbb{R}^{nm}$ and $\vec{a}=\text{vec}(\mat{A})=[{A}_{1,1}, {A}_{1,2}, \cdots, {A}_{m,n-1}, {A}_{m,n}]^T$, i.e., the concatenated vector of the transposed row vectors of $\mat{A}$.

\subsection{Optimal Transport and Unbalanced Optimal Transport}
\label{Sec:2.1}
{\bf Optimal Transport (OT):} Given two discrete probability measures $\vec{a}\in \R^{m}$ and $\vec{b} \in \R^{n}$ that $\|\vec{a}\|_1 = \|\vec{b}\|_1$, the standard OT problem seeks a corresponding {\it transport matrix} $\mat{T} \in \R_{+}^{m \times n}$ minimizing the total transport cost \cite{Kantorovich_1942}. This can be formulated as
\begin{eqnarray}
\label{Eq:Standard_OT}
\operatorname{OT}(\vec{a},\vec{b}) := \min_{ \mat{T} \in \R_{+}^{n \times m}} \langle \mat{C}, \mat{T} \rangle, \quad
\text{subject\ to} \quad \mat{T} \one_m= \vec{a}, \mat{T}^{T}\one_n = \vec{b},
\end{eqnarray}
where $\mat{C} \in \mathbb{R_{+}}^{n \times m}$ is the {\it cost matrix}. The famous Wasserstein q-distance is obtained for  $C_{i,j} = \|a_i - b_j\|_q$ \cite{Villani_2008_OTBook}.

As $\vec{t}\!=\!\text{vec}({\mat{T}}) \in \mathbb{R}^{nm}$ and $\vec{c}\!=\!\text{vec}({\mat{C}}) \in \mathbb{R}^{nm}$, we reformulate Eq.~(\ref{Eq:Standard_OT}) in a vector format as \cite{Chapel_NeurIPS_2021}
\begin{eqnarray}
\label{Eq:Vector_OT}
\operatorname{OT}(\vec{a},\vec{b}) := \min_{\vec{t} \in \R_{+}^{nm}} \vec{c}^T\vec{t}, \quad
\text{subject\ to}\quad \mat{N}\vec{t} = \vec{a}, \mat{M}\vec{t} = \vec{b}, \notag
\end{eqnarray}
where $\mat{N} \in \R^{n \times nm}$ and $\mat{M} \in \R^{m \times nm}$ are two {\it indices matrices} composed of ``0" and ``1" to compute row-sum and column-sum of $\mat T$. Each column of $\mat{N}$ and $\mat{M}$ has only a single non-zero element equal to $1$. Examples are listed in the supplementary. We denote $\hat{\vec t}$ as the {\it optimal primal solution} such that $\operatorname{OT}(\vec{a},\vec{b}) = \vec c^{T}\hat{\vec t}$.

{\bf Unbalanced Optimal Transport (UOT):}  The UOT {\it relaxes} the marginal constraints in Eq.~(\ref{Eq:Standard_OT}) by replacing the equality constraints with penalty functions on the marginals with divergence \cite{Caffarelli_AM_2010,chizat2017scaling}. Formally, defining $\vec{y} = [\vec{a}^T, \vec{b}^T]^T \in \mathbb{R}^{n+m}$ and the concatenation {\it index matrix} $\mat{X} = [\mat{M}^T,\mat{N}^T]^T \in \mathbb{R}^{(n+m) \times nm}$, the UOT can be formulated by introducing a penalty function for the discrete distributions as \cite{Chapel_NeurIPS_2021}
\begin{equation}
\label{eq:uot}
\operatorname{UOT}(\vec{a},\vec{b}) := \min_{\vec{t} \in \R_{+}^{nm}} \lambda\vec{c}^T\vec{t} + D_\phi(\mat{X}\vec{t},\vec{y}).
\end{equation}

\subsection{Duality and Strong Concavity}
We consider an optimization problem as
\begin{equation}
\label{eq:lassolike}
\min_{\vec{t} \in \mathbbm{R}^{n}} \left\{ f(\vec{t}) := g(\vec{t}) + h(\mat{X} \vec{t}) \right\},
\end{equation}
where $g: \mathbbm{R}^{n} \rightarrow (-\infty, \infty]$ and $h: \mathbbm{R}^{m} \rightarrow (-\infty, \infty]$ represent proper convex functions, $\vec{t} \in \mathbbm{R}^{n}$ is the primal optimization variable, and $\mat{X}\in \mathbbm{R}^{m\times n}$.
To derive the dual problem of Eq.~(\ref{eq:lassolike}), we rely on the Fenchel--Rockafellar Duality:

\begin{thm}[Fenchel--Rockafellar Duality {\cite{Rockafellar_Springer_1998}}]
\label{Thm:FRD}
Consider the problem Eq.~(\ref{eq:lassolike}). Assuming that all the assumptions are satisfied, we have
\begin{equation}
\label{Eq:FRD}
\min_{\vec{t}} g(\vec{t}) + h(\mat{X}\vec{t}) = \max_{\vec{\vec{\theta}}} -h^*(-\vec{\theta})-g^*(\mat{X}^T\vec{\theta}),
\end{equation}
where $\vec{\theta} \in \mathbbm{R}^m$, and where $g^*: \mathbbm{R}^n \rightarrow [-\infty, \infty]$ and $h^*: \mathbbm{R}^m \rightarrow [-\infty, \infty]$ respectively stand for the convex conjugate of the extended real-valued functions $g: \mathbbm{R}^{n} \rightarrow (-\infty, \infty]$ and $h: \mathbbm{R}^m \rightarrow (-\infty, \infty]$. We call the former primal problem $\mat P(\vec t)$ and the latter dual problem $\mat D(\vec \theta)$.
\end{thm}

Next, we introduce the concepts of strong concavity and $\mat L$-blockwise strong concavity, the latter being a more generalized form.
\begin{definition}[Strong Concavity]
\label{Def:SC}
Function $f(\vec x)$ is an {\it L-strongly concave} function if there exists a positive constant L that, for $\forall \vec y \in \mathbbm{R}^{n}$:
\begin{equation}
\label{eq:sc}
f(\vec y) \leq f(\vec x) +\nabla f(\vec x)^{T}(\vec y-\vec x) - \frac{L}{2}\|\vec y -\vec x\|_2^2.
\end{equation}
\end{definition}

\begin{definition}[$\mat L$-Blockwise Strong Concavity]
\label{Def:BSC}
A function $f(\vec x)$ is $\mat L$-blockwise strongly concave if there exists a diagonal positive-definite matrix $\mat L$ such that
\begin{equation}
\label{eq:bsc}
f(\vec y) \leq f(\vec x) +\nabla f(\vec x)^{T}(\vec y-\vec x) - \frac{1}{2}(\vec y -\vec x)^{T}\mat L (\vec y -\vec x).
\end{equation}
\end{definition}

The $\mat L$ of $f(\vec x)$ can be obtained by computing a lower bound on the Hessian matrix of the $f(\vec x)$ \cite{DBLP:books/cu/BV2014,DBLP:journals/siamjo/BeckT13}. Assuming that $\mat L$ exhibits no anisotropy. In this case, $\mat L = L\textbf{I}$. When Eq.~(\ref{eq:sc}) and Eq.~(\ref{eq:bsc}) only hold for $\vec y \in \mathcal{G} \subset \mathbbm{R}^{n}$, they referred to as {\it $\mathcal{G}$-locally L-strongly concave} and {\it $\mathcal{G}$-locally $\mat L$-blockwise strongly concave}, respectively.

\section{UOT SAFE SCREENING}
\label{sec:pro}

In high-dimensional regression tasks, regularization terms like $\ell_1$-norm bring desirable sparsity to solutions. {\it Safe Screening} seeks to extend this benefit by identifying and discarding zero values in the solution before the computationally expensive optimization process concludes. It was applied to accelerate Lasso problems \cite{ghaoui2010safe} and SVMs \cite{Ogawa_ICML_2013} by eliminating unnecessary features. In recent years, researchers have continuously proposed new methods to narrow down the safe region while ensuring acceptable computational costs \cite{Liu_ICML_2014,Wang_JMLR_2015, JMLR:v18:16-577, Yamada_NIPS_2021}. Furthermore, they actively extend the application of Safe Screening to a broader spectrum of functions and regularizations. \cite{DBLP:conf/icml/AtamturkG20, Dantas_ICASSP_2021, 10.5555/3546258.3546494}

This section presents the Safe Screening framework proposed for the UOT. We begin by formulating the dual form $\mat D(\vec \theta)$ and discussing the {\it dual feasible region}, denoted as $\mathcal{R}^D$, we then propose {\it Shifting Projection} based on the special index matrix $\mat X$. We discuss the strong concavity under different penalty functions and explain how they influence the design of the safe region $\mathcal{R}^{S}$. Next, we focus on the KL-penalty function, proving its {\it $\mat L$-blockwise local strong concavity} and introducing the {\it Gap-Ellipse} (Ell) method. Finally, considering the unique structural characteristic of $\mat X$ in the UOT, we design and propose the Cruciform two-hyperplane (CTP) method to further enhance the effectiveness of Safe Screening. Concrete proofs are provided in the supplementary material.

\subsection{Dual Formulation and Feasible Region}
\begin{lem}[Dual Form of UOT]
\label{Lem:alluotdual}
Applying {\bf Theorem~\ref{Thm:FRD}}, we can derive the dual problem for UOT as follows:
\begin{eqnarray}
\max_{\vec{\theta} \in \R^{n+m}} &\!\!\!&\mat D(\vec \theta) = -\displaystyle{D_{\phi}^*(-\vec{\theta})}  \nonumber \\ 
\text{subject\ to}&\!\!\!&\vec{x}_p^T\vec{\theta} -\lambda c_p \leq 0, \quad \forall p \in [nm],
\label{eq:alluotdual}
\end{eqnarray}
\end{lem}

\begin{table}[t]
\caption{Characteristics of dual function and Hessian matrix in UOT under different penalty Functions.}
\label{table:sc}
\begin{center}
\begin{tabularx}{\textwidth}{@{}lXlll@{}}
\toprule
Penalty & \!\!$D_\phi(\mat{X}\vec{t},\vec{y})$ & $\mat D(\vec \theta)$ & \!\! ($\partial^2\! \mat D\!/\partial \vec\theta^2$) & \!\!Dual Characteristic \\ \hline  \addlinespace[3pt]
$\ell_2$ & \!\!$\|\mat X\vec t - y\|_2^2$ & $-\frac{1}{2}\|\vec{\theta}\|_2^2+\vec{y}^T\vec{\theta}$ & \!\!\!\!$-\operatorname{Diag}(\one)$ & \!\!Strongly concave \\  \addlinespace[3pt]\hline \addlinespace[5pt]
KL & \makecell[l]{$\!\!(\mat X\vec t + \epsilon )\ln ((\mat X \vec t \!+\!\epsilon)/ \vec y)$ \\ \!\!$-(\mat X\vec t +\epsilon) + \vec y$} & \makecell[l]{$-\vec y^{T} e^{ - \vec \theta} + \vec y^{T} \one_{n+m} $\\$-\epsilon\vec\theta^{T}\one_{n+m}$} & \!\!\!\!$-\operatorname{Diag}(\vec y^{T}\!\!\odot\! e^{-\vec\theta})$ & \!\!Not strongly concave \\ \addlinespace[3pt] \hline \addlinespace[3pt]
TV & \!\!$\|\mathbf{X}\vec{t} - \vec{y}\|_1$ & \multicolumn{1}{@{}c@{}}{
  $\begin{array}{@{}l@{}}
    \begin{cases}
      \vec{\theta}^{T}\vec{y}, \forall p,~ |\theta_p| < 1, \\
      +\infty, \text{ otherwise}.
    \end{cases}
  \end{array}$} & \!\!\!\!$\mathbf{O}$ & \!\!Linear \\ \addlinespace[1pt]
\bottomrule
\end{tabularx}
\end{center}
\end{table}

Different choices of penalty functions result in different dual functions (see Table~\ref{table:sc}), but the dual feasible region, denoted as $\mathcal{R}^{D}$, composed of the constraints in Eq.~(\ref{eq:alluotdual}), remains consistent, and the { \it optimal dual solution} $\hat{\vec{\theta}} \in \mathcal{R}^{D}$.

According to the KKT conditions of the Fenchel--Rockafellar duality in {\bf Theorem~\ref{Thm:FRD}}, a connection exists between the $\hat{\vec{t}}$ and $\hat{\vec{\theta}}$ as presented below.
\begin{thm} (Primal-Dual Relationship of Optimal Solution in the UOT) For the optimal primal solution $\hat{\vec{t}}$ and the optimal dual solution $\hat{\vec{\theta}}$, we have the following relationship for $\forall p, ~p \in [nm]$:
\label{Thm:KKT}
\begin{equation}
\label{eq:kkt}
\begin{split}
\vec{x}_p^T\hat{\vec{\theta}} -\lambda c_p \left\{
\begin{aligned}
< 0, \quad& \Longrightarrow \quad\hat{t}_p = 0,\\
= 0, \quad& \Longrightarrow \quad\hat{t}_p \geq 0.
\end{aligned}
\right.
\end{split}
\end{equation}
\end{thm}

It is important to note that $\vec{x}_p^T\hat{\vec{\theta}}$ on the left-hand side (LHS) of Eq.~(\ref{eq:kkt}) differs from $|\vec{x}_p^T\hat{\vec{\theta}}|$ in the standard Lasso problem due to the constraint $t_p > 0$ in the UOT problem. Hence, the fundamental concept behind the Safe Screening is to identify a region $\mathcal{R}^{S}$ such that $\hat{\vec{\theta}} \in \mathcal{R}^S$. If the following inequality holds:
\begin{equation}
\label{eq:kktineq}
\max_{\vec{\theta} \in \mathcal{R}^S} \vec{x}_p^T\vec{\theta} -\lambda c_p < 0,
\end{equation}
then we have $\vec{x}_p^T\hat{\vec{\theta}} -\lambda c_p < 0$, and the corresponding $\hat{{t}}_p = 0$. Consequently, $t_p$ can be screened out. Safe Screening seeks a {\it safe region} $\mathcal{R}^{S}$ containing $\hat{\vec{\theta}}$ as small as possible to make accurate and anticipatory judgments. The tool we rely on to construct this is the strong convexity in {\bf Definition~ \ref{Def:SC}}.

\subsection{Shifting Projection}
\label{Sec:SP}
To construct $\mathcal{R}^S$, a {\it feasible dual variable}, denoted $\tilde{\vec{\theta}}$, belonging to the set $\mathcal{R}^{D}$ is required. The variable $\vec \theta$, derived from $\vec t$ through the Fenchel--Rockafellar relationship in Eq.~(\ref{Lem:alluotdual}), may not satisfy the dual constraints, resulting in $\vec{\theta}\notin\mathcal{R}^{D}$ and necessitating projections. As illustrated in {\bf Algorithm~ \ref{Alg:UOTDynamicScreening}}, as $\vec t^{k}$ approaches $\hat{\vec t}$ during optimization, $\vec \theta^{k}$ similarly converges to $\hat{\vec \theta}$. A poor projection may significantly deviate from $\hat{\vec{\theta}}$, leading to an expansion of $\mathcal{R}^{S}$ and adversely affecting the screening outcome. This necessitates a high degree of accuracy in the projection. However, since $\mathcal{R}^{D}$ in Eq.~(\ref{eq:alluotdual}) is a polytope composed of $nm$ hyperplanes, regular projections can be computationally expensive. To address this issue, researchers have employed {\it Residuals Rescaling} \cite{DBLP:conf/icml/MassiasSG18} to achieve cost-effective approximate projections \cite[Proposition 11]{JMLR:v18:16-577} \cite[Theorem 11]{Yamada_NIPS_2021}.

\begin{figure}[t]
\centering
\begin{minipage}[t]{0.48\textwidth}
  \centering
  \includegraphics[width=\linewidth]{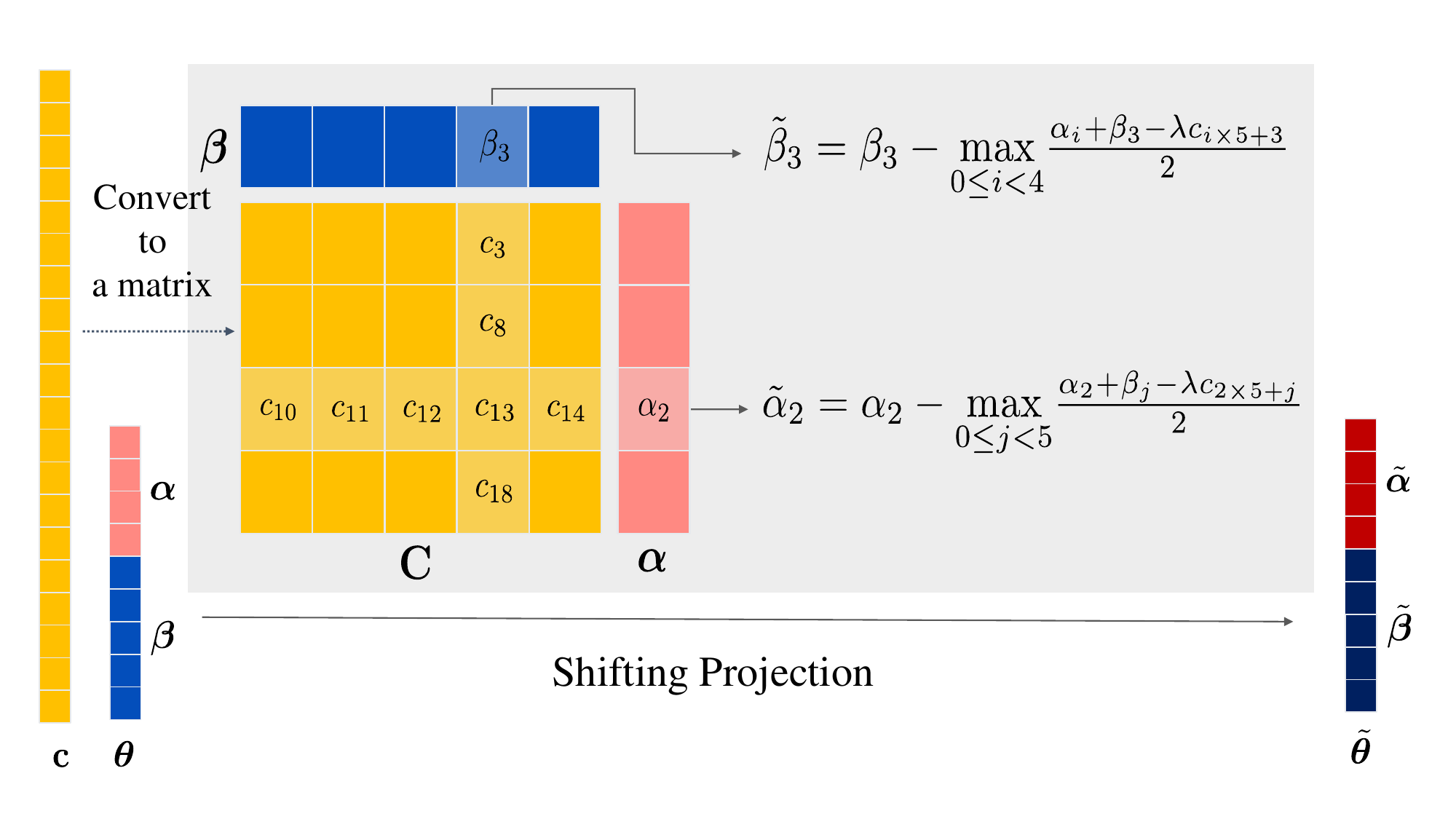}
  \caption{Example of Shifting Projection. Arranging $mn$ constraints into a yellow matrix. The light yellow blocks show the constraints that directly affect the corresponding elements.}
  \label{Fig:structure}
\end{minipage}
\hfill
\begin{minipage}[t]{0.48\textwidth}
  \centering
  \includegraphics[width=1.0\linewidth]{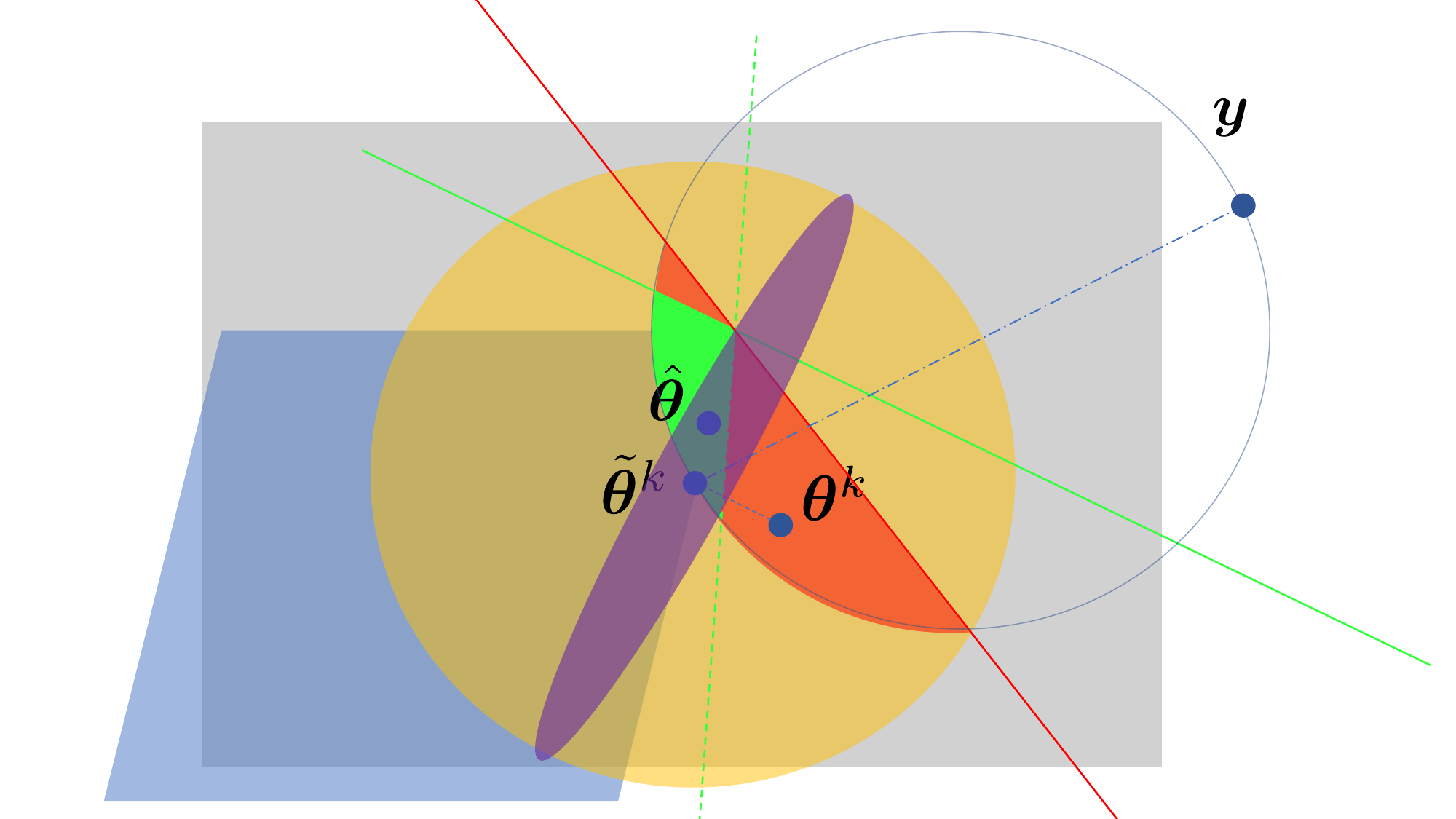}
  \caption{Comparison of different regions. $\mathcal{R}^{D}$(blue), $\mathcal{R}^{B}$(grey), Gap (yellow), Sa (orange), Sa-CTP (green), and Gap-Ell (purple).}
  \label{Fig:area}
\end{minipage}
\end{figure}

One extension of it to UOT can be achieved by considering $\tilde{\vec{\theta}} = \vec\theta /\max(1, \| \mat{X}^T\vec\theta/\lambda \vec{c}\|_{\infty})$, where the division $(\cdot\ /\cdot)$ is performed element-wise. The computational complexity is $O(\max(n,m)(n+m))$. However, Residuals Rescaling in UOT can exhibit significant fluctuations due to the varying range of values $c_p$ in the denominator term, and it may degenerate when there exists a $p$ such that $c_p = 0$. To overcome this issue, we propose the {\it Shifting Projection}:

\begin{thm}[Shifting Projection for the UOT]
\label{Thm:UOT_ShiftProjection}
Let $\vec{\theta} = [{\vec{\alpha}}^T,{\vec{\beta}}^T]^T$, where $\vec{\alpha}\in\R^{m}$ and $\vec{\beta}\in\R^{n}$. For $u \in [m]$ and $v \in [n]$, the Shifting Projection is defined as:
\begin{equation}
\left\{
\begin{array}{ll}
\label{eq:uotproj}
\tilde{{\alpha}}_u &= \displaystyle{{{\alpha}}_u - \max_{0\leq j < n} \frac{{{\alpha}}_u +{{\beta}}_j - \lambda{c}_{un+j}}{2}} \\
\tilde{{\beta}}_v &= \displaystyle{{{\beta}}_v - \max_{0 \leq i < m} \frac{{{\alpha}}_i +{{\beta}}_v - \lambda{c}_{in+v}}{2}}.
\end{array}
\right.
\end{equation}
\end{thm}
Since $\vec x_p$ is an index vector containing only two elements with a value of "1", we can express the constraint as $\vec x_p^{T}\vec\theta = \alpha_{u} + \beta_{v} < \lambda c_p$, where $(u,v)=(p \mid n,p \mod n)$. In Figure \ref{Fig:structure}, it can be observed that the variables need to be shifted by half of the maximum value among the row and column constraints. The Shifting Projection is particularly suitable for UOT as it avoids degeneracy issues and takes into account differences in coordinate dimensions.

\subsection{Safe Region Construction}
\label{Sec:Ell}

After obtaining $\tilde{\vec{\theta}}\in\mathcal{R}^{D}$, we construct $\mathcal{R}^{S}$ to include $\hat{\vec{\theta}}$ if the dual function satisfies {\bf Definition~\ref{Def:SC}}. We choose simple geometries, such as balls and domes for testing Eq.~(\ref{eq:kktineq}) due to their computational efficiency requirements.
\label{Sec:SC}
\begin{thm}[Ball Region Construction \cite{Yamada_NIPS_2021}]
\label{Thm:Gap_all}
For a $\mathcal{R}^{D}$-locally $L$-strongly concave function $\mat D(\vec\theta)$, finding a point $\tilde{\vec\theta} \in \mathcal{R^{D}}$, then we can construct the following region that guarantees $\hat{\vec\theta} \in\mathcal{R}^{G}$:
\begin{equation}
\mathcal{R}^{G}(\tilde{\vec\theta}) :=\left\{\vec\theta \ |\|(\vec\theta-\tilde{\vec\theta})\|_2^2 \leq \frac{2(\mat D(\vec \theta) - \mat D(\tilde{\vec\theta}))}{ L} \right\}.
\label{eq:gs}
\end{equation}
\end{thm}
The famous Gap-safe method \cite{JMLR:v18:16-577} uses a constant bound to relax the right-hand side (RHS): $\mat D(\vec \theta) - \mat D(\tilde{\vec\theta}) < \mat \mat P({\vec t}) - D(\tilde{\vec \theta}) = \Gap(\vec t, \tilde{\vec \theta})$, which is known as {\it weak duality}. Sasvi method in \cite[Theorem 8]{Yamada_NIPS_2021} computes $\mathcal{R}^{S}(\tilde{\vec\theta}) :=\{\vec\theta \ | (\vec\theta-\tilde{\vec\theta})^{T}(\vec\theta-\vec y)\leq 0\}$, which is proved in Theorem 6 of \cite{Yamada_NIPS_2021}, are both {\it balls}.

Considering the anisotropy at different coordinates of the dual KL-penalized UOT problem, we propose a new screening method called {\it Gap-Ellipse Safe Screening} method:
\begin{thm}[Gap-Ellipse Safe Screening]
\label{Thm:Gap_ellipse}
For a $\mat L$-blockwise strongly concave function $\mat D(\theta)$, finding a point $\tilde{\vec\theta} \in \mathcal{R^{D}}$, then we can construct the following ellipse region that guarantees $\hat{\vec\theta} \in\mathcal{R}^{G}$:
\begin{equation}
\mathcal{R}^{G}(\tilde{\vec\theta}) :=\{\vec\theta \ | (\vec\theta-\tilde{\vec\theta})^{T}\mat L(\vec\theta-\tilde{\vec\theta}) \leq 2\Gap(\vec t, \tilde{\vec \theta})\}.
\label{eq:ge}
\end{equation}
\end{thm}
Considering the dual function of the UOT, shown in Table~\ref{table:sc}, the TV-penalty results in linearity, while the L2-penalty leads to strong convexity. However, the KL-penalty does not exhibit strong convexity. Only the $\ell_2$-penalty can apply the {\bf Theorem~\ref{Thm:Gap_all}}. Therefore, we consider proving the local strong concavity of the dual function in a new region under the conditions of the KL-penalty.

If the dual function exhibits local strong concavity on $\mathcal{R}^{D}$, then {\bf Theorem~\ref{Thm:Gap_all}} remains applicable. The dual function is to be $\mathcal{R}^{D}$-locally strongly concave when the index matrix $\mat X$ has a right inverse \cite{Dantas_ICASSP_2021}. However, in the case of UOT, the index matrix $\mat X$ does not possess a right inverse, as proven in the supplementary material. This prompts us to seek an alternative approach to establish its locally strong concave characteristic by leveraging the unique structure of UOT.

We propose to construct a new {\it optimal feasible region}, $\mathcal{R}^{B}$, such that $ \hat{\vec\theta} \in \mathcal{R}^{B}$. We then demonstrate that the KL-penalized UOT has a locally strongly concave dual function on $\mathcal{R}^{B} \cap \mathcal{R}^{D}$.

\begin{thm}[Bounds of $\hat{\vec \theta}$ for KL-penalized UOT]
\label{Thm:kl_local}
For $\hat{\vec{\theta}} = [\hat{\vec{\alpha}}^T,\hat{\vec{\beta}}^T]^T$ and $(u,v)=(p \mid n,p \mod n)$ as stated in {\bf Theorem~\ref{Thm:UOT_ShiftProjection}}, considering the symmetry of the rows and columns, we focus on $\alpha_u$ and use $\operatorname{Low}(\vec \theta,v+n)$ to represent a lower bound of $\beta_v$ derived from $\forall \vec{\theta}$, without loss of generality. Thus, we have
\begin{equation*}
\begin{split}
\hat{\alpha_u} < \overline{\alpha_u}
< \min_v (\lambda c_p -\underline{\beta_v} ).
\end{split}
\end{equation*}
Here, we use $\overline{\alpha_u}$ to represent its upper bound, and $\underline{\beta_v} := \operatorname{Low}(\vec \theta, j)$, with $j = n+v$, to represent the lower bound. We define the lower bound obtained from $\vec\theta$ as:
\begin{equation*}
\begin{split}
\operatorname{Low}(\vec \theta, j) = \ln\left(\frac{\epsilon - y_j}{K^j - y_j +\epsilon}\right),
\end{split}
\end{equation*}
where $K^{j}= \mat D(\vec\theta) - \mat P^{j}({\vec t})$, and $\mat P^{j}$ is the primal problem of $\sum_{i\neq j}^{n+m} \mat D^{i}(\theta_i)$. The dual problem has a decomposable structure where $\mat D(\vec{\theta})=\sum_{i=1}^{n+m}\mat D^i(\theta_i)$, and $\mat D^i(\theta_i)$ represents the value of the $i$-th coordinate of the dual problem.
\end{thm}

We construct a computable primal function $\mat P^{j}({\vec t})$ and utilize it to provide new bounds $\underline{\vec\theta}$ and $\overline{\vec\theta}$ for $\hat{\vec\theta}$ as $\vec\theta$ updates during optimization. Detailed proof is provided in the supplementary material. By using these bounds, we construct a Box region $\mathcal{R}^{B}(\vec\theta,\vec\theta^{k}) = {\vec\theta | \underline{\theta_i} \leq {\theta}_i \leq \overline{\theta}_i, \ i \in [n+m] }$, and ensure that $\hat{\vec\theta} \in \mathcal{R}^{B}$. We can find $\tilde{\vec\theta} \in \mathcal{R}^{B}\cap \mathcal{R}^{D}$ as projecting on Box region is a coordinate-wise shifting, given that the Hessian function is an increasing function, we can demonstrate that $\mat L = \diag{(\vec y e^{-\overline{\vec\theta}})}$, and the dual function of the KL-penalized UOT problem is ($\mathcal{R}^{B}\cap \mathcal{R}^{D}$)-locally $\mat L$-blockwise strongly concave.

The LHS in Eq.~(\ref{eq:ge}) represents an ellipse rather than a ball with the center of $\tilde{\vec\theta}$. {\bf Theorem~\ref{Thm:kl_local}} enables us to take into account the anisotropic property of the function as diagnose elements in Hessian hugely vary, which leads to considerable performance improvement in the screening process.

\paragraph{Two-hyperplane Safe Screening Region.}
\label{Sec:CTP}

We can adapt the ball or ellipse $\mathcal{R}^{G}$ as $\mathcal{R}^{S}$, and combine it with the $\mathcal{R}^{D}$ to construct a smaller $\mathcal{R}^{S}$. As $\mathcal{R}^{D}$ is a $(n+m)$-polytope with $nm$ facets, finding the maximum value within it is expensive \cite{Wang_JMLR_2015}. A simple approach is to relax all the hyperplanes into a single hyperplane, transforming the $\mathcal{R}^{S}$ into the intersection of a sphere or ellipse with the half-space \cite{DBLP:conf/icassp/XiangR12, DBLP:conf/icml/LiuZWY14, Yamada_NIPS_2021, JMLR:v18:16-577}. This allows us to obtain the maximum value at an acceptable cost. We denote the relaxed region by $\mathcal{R}^{Re}$. and then the new region can be constructed by $\mathcal{R}^{S} = \mathcal{R}^{G} \cap \mathcal{R}^{Re}$.

\begin{thm}[Dome safe region for UOT \citep{Yamada_NIPS_2021}]
\label{Thm:Dome}
For every primal variable $t_p \geq 0$, $p \in [mn]$, we construct a specific region $\mathcal{R}^{Re}_{p}$ as presented below.
\begin{equation}
\label{Eq:FinalDome}
\begin{split}
\mathcal{R}^{Re}(\vec{\theta}, \vec{t}):= \left\{\vec{\theta} \ \left|\
\begin{aligned}
\sum_{p\in [nm]}(\vec{x}_{p}^{T}\vec{\theta} - \lambda {c}_p)t_p\leq 0 \\
\end{aligned}
\right.
\right\}.
\end{split}
\end{equation}
\end{thm}

Considering the specific structure of $\mat X$ in the UOT, we propose to relax the polytope into the intersection of half-spaces by combining the $nm$ dual constraints with a positive weight $t_p \geq 0$.

\begin{thm}[Cruciform two-hyperplane (CTP) safe region for UOT]
\label{Thm:AreaScreeninguOT}
For every primal variable $t_p \geq 0$, $p \in [mn]$, let $I_p = \{ i \ | \ 0\leq i<nm, u = i\mid n \text{ and } v = i\mod n\}$, and ${I}^{C}_p = \{ i \ | \ 0\leq i<nm, i \notin I_p\}$. Then we construct a specific region $\mathcal{R}^{Re}_{p}$ as presented below.
\begin{equation}
\label{Eq:FinalRS}
\begin{split}
\mathcal{R}^{Re}_{p}(\vec{\theta}, \vec{t}):= \left\{\vec{\theta} \ \left|\
\begin{aligned}
\sum_{p\in I_p}(\vec{x}_{p}^{T}\vec{\theta} - \lambda {c}_p)t_p\leq 0, \sum_{p\in {I}^{C}_p}(\vec{x}_{p} ^{T}\vec{\theta}- \lambda {c}_p)t_p\leq 0 \\
\end{aligned}
\right.
\right\}.
\end{split}
\end{equation}
\end{thm}

To elaborate, we partition the polytope constraints into a primary group $I_p$ and a secondary group ${I}^{C}_p$. For each $p$, $\mathcal{R}^{Re}_{p}$ is the intersection of two relaxed half-spaces. The effectiveness of $\mathcal{R}^{Re}_{p}$ can be intuitively understood: while more hyperplanes result in a tighter region, in the standard Lasso problem, given that all optimization directions in $\vec x_p^{T} \theta$ are of equal importance, an arbitrary split is likely to result in nearly {\it parallel} high-dimensional hyperplanes. This phenomenon is illustrated in Figure~\ref{Fig:ex2} in {\bf Section~\ref{Sec:EX2}}, where using random hyperplanes does not lead to improved screening.

In contrast, we can incorporate $n+m -1$ constraints, depicted in light yellow in Figure~\ref{Fig:structure}, to form the primary hyperplane. These constraints are directly associated with the optimization direction $\alpha_u$ and $\beta_v$. The remaining constraints are grouped as the secondary hyperplane, whose coefficients of $\alpha_u$ and $\beta_v$ would be zero. This arrangement avoids parallelism with the primary hyperplane and significantly shrinks the safe region in comparison with the Dome method in {\bf Theorem~\ref{Thm:Dome}}.

It is worth noting that, unlike \cite{7469344}, we design two specific hyperplanes for each element in the UOT and maximize Eq.~(\ref{eq:kktineq}) on them. Consequently, as in Figure~\ref{Fig:area}, the two green lines represent two hyperplanes. Optimizing on $\mathcal{R}^{Re}_{p} \cap \mathcal{R}^{G}$ has a closed-form solution by using the Lagrangian method. Its computational process is in the supplementary material.

\subsection{Optimization Algorithm and Computational Cost Analysis}

\begin{figure}[t]
\begin{minipage}[t]{0.55\textwidth}
\begin{algorithm}[H]
\caption{UOT with Safe Screening}
\begin{algorithmic}
\label{Alg:UOTDynamicScreening}
\renewcommand{\algorithmicrequire}{\textbf{Input:}}
\renewcommand{\algorithmicensure}{\textbf{Output:}}
\REQUIRE $\vec{t}^0, \vec c, \vec s \in \mathbb{R}^{nm}, \mat X \in{\mathbbm{R}^{(n+m)\times nm}}, s_{p}=1\ ,\forall p \in [mn], w \in \mathbbm{N}$
\ENSURE $\vec t^{K}$
\FOR {$k = 0 \text{ to } K-1$}
\IF{($k\mod w) =0$}
\STATE $\tilde{\vec{\theta}}^{k} = \operatorname{Projection}(\vec{\theta}^k)$ \hfill Eq.~(\ref{eq:uotproj})
\FOR {$p = 0 \text{ to } \text{length}(\vec t)$}
\STATE $\mathcal{R}^{S} \leftarrow \mathcal{R}^{G}(\tilde{\vec \theta^{k}})\cap \mathcal{R}^{Re}{(\vec{\theta}^{k},\vec{t}^k)}$ \hfill Eq.~(\ref{Eq:FinalDome},\ref{Eq:FinalRS})
\STATE
\STATE{$\vec s \leftarrow {s_{p} = 0 \text{ if } {\displaystyle \max_{\vec{\theta} \in \mathcal{R}^S} {\vec x_{p}}^T\vec{\theta}^{k} <\lambda c_{p} }}$ \hfill Eq.~(\ref{eq:kktineq})}
\STATE \text{Delect }$ t_p^{k}, c_p^{k} \text{ and } \vec x_p \text{ if } s_{p}\neq 0$
\ENDFOR
\ENDIF
\STATE $\vec{t}^{k+1} = \operatorname{update}(\vec{t}^k, \vec c, \mat X)$
\STATE $\vec s = \one_{\text{length}(\vec t)}$
\ENDFOR
\end{algorithmic}
\end{algorithm}
\end{minipage}
\hfill
\begin{minipage}[t]{0.4\textwidth}
\label{Fig:ex2_1}
\centering
\includegraphics[width = 1\linewidth]{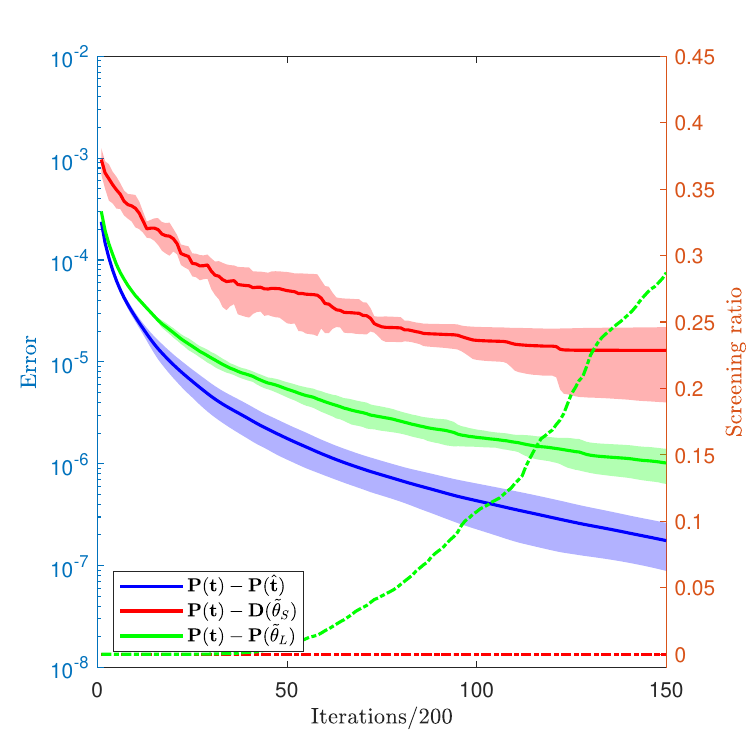}
\caption{Comparison of the dual gap (shadows) between Shifting Projection (Green) and Residual Rescaling (Red) with the Objective Error (Blue), and its influence on Screening Ratio (lines).}
\label{Fig:ex1}
\end{minipage}
\end{figure}

\paragraph{Safe Screening Algorithm.} We define the {\it screening mask vector}, denoted as $\vec{s} \in \mathbb{R}^{mn}$ at first, for the transport vector $\vec{t}$. We remove $t_p$ if $s_p = 0$ before sending it to the optimizer. The entire optimization process is summarized as in {\bf Algorithm \ref{Alg:UOTDynamicScreening}}. The $\operatorname{update}(\vec{t},\vec c, \mat X)$ operator in the algorithm indicates that the optimizer updates $\vec{t}^k$ into $\vec{t}^{k+1}$. The dimension of $\vec t^{k}$ would gradually decrease as well as $\vec s$. $w \in \mathbb{N}$ specifies the screening period to update the mask $\vec s$ according to the intensity of the change on the screening ratio, decided by the user. Faster convergence of the optimizer is associated with a smaller period $w$. It must be emphasized that the Safe Screening is {\it independent} of the adopted optimization algorithm.

\paragraph{Computational Cost Analysis.}
The Shifting Projection has the same computational burden as Residuals Rescaling for the UOT. The Gap-Ellipse method has the same expense as Gap as $O(nm)$ because only two dual coordinates are activated each time. The CTP method must apply the Lagrangian method for every primal element $t_p$. However, under the special structure of the $\mat X$, for every $t_p$, the data required during the Lagrangian method can be divided as row summation and column summation of the transport matrix $\mat T$, which can be computed together and be reused for other elements that have the same $t \mod m $ or $t \mid m$. This process of summarization helps to preserve the overall optimization complexity to $O(mn)$. The detailed steps are in the supplementary material.

\section{EXPERIMENTS}

\begin{figure*}[t]
\centering
\includegraphics[width = \textwidth]{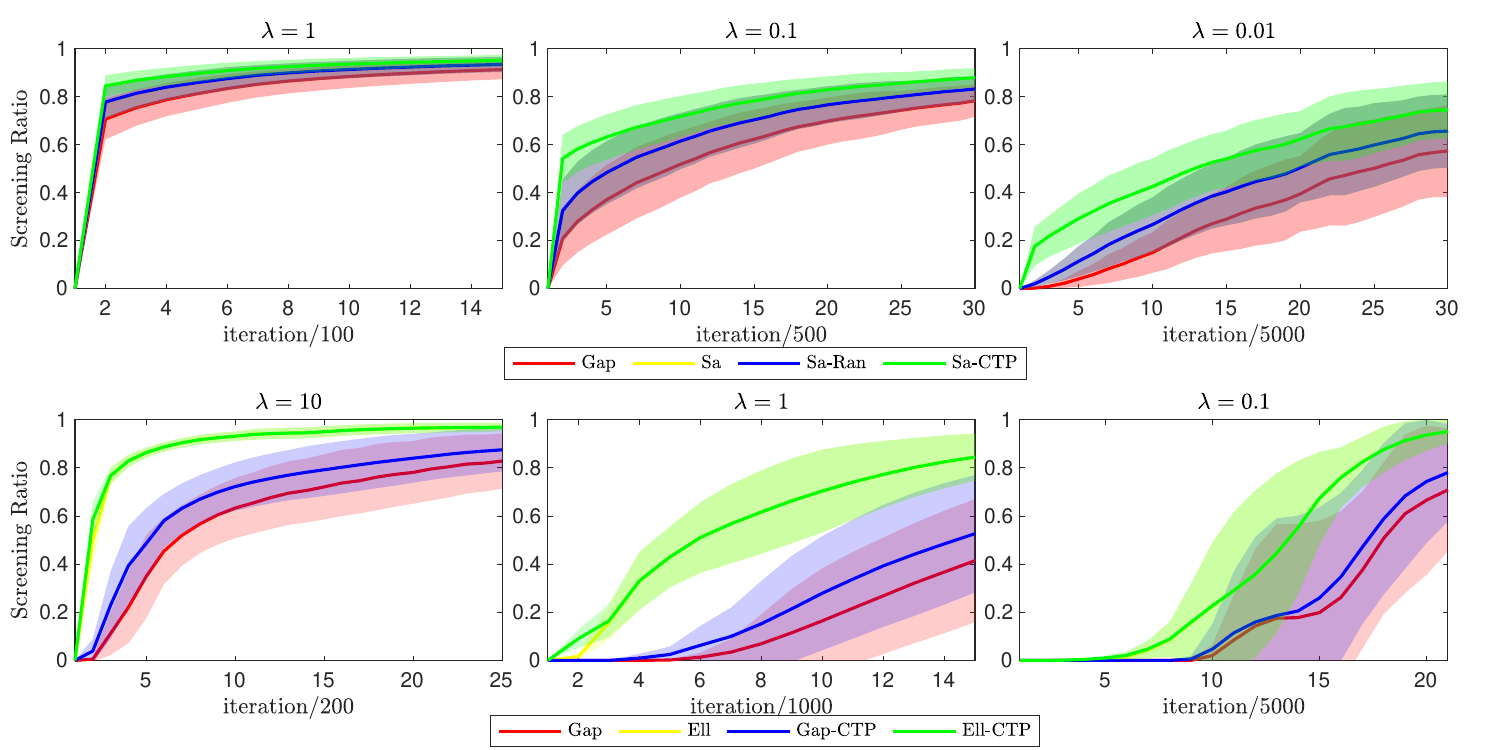}
\caption{The Screening for different $\lambda$ in $\ell_2$ (Upper) and KL-penalized (Lower) UOT.}
\label{Fig:ex2}
\end{figure*}

\begin{table}[t]
\caption{Speed-up ratios of the Safe Screening: FISTA for $\ell_2$ penalty, and MM for KL penalty. The number in the second row represents dual gap $\varepsilon$.}
\label{sample-table}
\centering
\begin{tabular}{l|c|cc|cc|c|cc|cc}
\toprule
Dataset & Method & \multicolumn{4}{c|}{$\ell_2$} & Method & \multicolumn{4}{c}{KL} \\
& & \multicolumn{2}{c}{$\lambda=10^{-1}$} & \multicolumn{2}{c|}{$\lambda=10^{-2}$} & & \multicolumn{2}{c}{$\lambda=1$} & \multicolumn{2}{c}{$\lambda=10^{-1}$} \\
\cmidrule{3-6} \cmidrule{8-11}
& & $10^{-5}$ & $10^{-7}$ & $10^{-5}$ & $10^{-7}$ & & $ 10^{-5}$ & $10^{-7}$ & $ 10^{-5}$ & $10^{-7}$ \\
\midrule
MNIST & Gap & 1.55 & 3.11 & 1.01 & 1.31 & Gap & 0.97 & 1.18 & 0.99 & 1.00 \\
& Sa & 1.78 & 3.73 & 1.02 & 1.48 & Gap-CTP& 0.97 & 1.28 & $\bm{1.01}$ & 1.01 \\
& \!Sa-CTP\! & $\bm{2.14}$ & $\bm{4.51}$ & $\bm{1.11}$ & $\bm{1.82}$ & Ell & $\bm{1.23}$ & $\bm{2.13}$ & 1.00 & $\bm {1.11}$ \\
& & & & & & \!\!Ell-CTP\!\! & $\bm{1.23}$ & 2.11 & 1.00 & 1.10 \\
\midrule
Gauss & Gap & $\bm {1.58}$ & $\bm {3.96}$ & 1.00 & 1.63 & Gap & 1.02 & 1.29 & 0.99 & 0.98 \\
& Sa & 1.41 & 2.78 & 1.01 & 2.19 & Gap-CTP & 1.00 & 1.49 & $\bm{1.00}$ & 0.99 \\
& Sa-CTP & 1.46 & 2.87 & $\bm{1.02}$ & $\bm{2.93}$ & Ell & 1.75 & $\bm{4.96}$ & 0.99 & $\bm{1.42}$ \\
& & & & & & \!\!Ell-CTP\!\! & $\bm{1.78}$ & 4.80 & 0.99 & 1.41 \\
\bottomrule
\end{tabular}
\end{table}

\label{sec:exp}
We conduct three experiments on randomly generated high dimensional Gaussian distributions and MNIST dataset. We first validate the effectiveness of the Shifting Projection in Eq.~(\ref{Thm:UOT_ShiftProjection}). Next, we compare the screening ratio of different methods for the $\ell_2$ and KL-penalized problems. Finally, we apply our Safe Screening on some popular optimization solvers, i.e., FISTA \cite{Beck_2009_SIAMIS}, Majorization-Minimization(MM) \cite{Chapel_NeurIPS_2021}, coordinate descent, and BFGS for the $\ell_2$-norm penalized UOT, and MM, and BFGS method for the KL-penalized problem to elucidate its computation speed-up ratios. All the code was written in Matlab. Due to space constraints, we have included the experiments involving coordinate descent and BFGS in the supplementary materials.

\paragraph{Projection Method Comparison.}
\label{Sec:EX1}
We evaluate the dual gap of two approximate projection methods as depicted in Figure \ref{Fig:ex1}. We utilize the FISTA algorithm to solve 10 pairs of $\ell_2$-penalized 100-dimensional Gaussian distribution transport problems. $\vec\theta_S$ and $\vec\theta_L$ are obtained by Shifting Projection and Residuals Rescaling with $\lambda = 10^{-2}$, respectively. We then calculate the corresponding dual gap. Our method is more accurate and stable. However, Residuals Rescaling exhibits slow screening progress due to the subpar dual gap.

\paragraph{Screening Method Comparison.}
\label{Sec:EX2}
We evaluate the screening ratios of different methods on the $\ell_2$-penalized and KL-penalized UOT. For the KL-penalized problem, we select $\lambda \in \{10^{-1}, 10^{0}, 10^{1}\}$, and for the $\ell_2$-penalized problem, we select $\lambda \in \{10^{0}, 10^{-1}, 10^{-2}\}$. All methods compared employ our Shifting Projection as it has been demonstrated to be more suitable for the UOT than other techniques in previous experiments. We conducted experiments on the MNIST dataset 5 times, solving them with the FISTA and MM algorithms. For the KL-penalized problem, we compare the Gap (Gap), Gap with CTP (Gap-CTP), Gap-Ellipse (Ell), and Gap-Ellipse with CTP (Ell-CTP). For the $\ell_2$-penalized problem, we select the Gap (Gap), Dynamic Sasvi-Dome (Sa), Sasvi with CTP (Sa-CTP), and Sasvi with randomly-selected two-hyperplanes (Sa-Ran) methods. The results are depicted in Figure \ref{Fig:ex2}. The blue shadow almost overlaps with the yellow Sa method. This suggests that the simple addition of more planes is not the key factor contributing to the success of the CTP method. Rather, the advantage arises from the specific composition of the planes.

\paragraph{Computational Speed-up Ratio.}
\label{para:ex3}
In this section, we conduct experiments on both the MNIST and Gaussian datasets using different screening methods 5 times to demonstrate the speed-up ratio of these methods on $\ell_2$ and KL-penalized UOT. We adopt two stopping conditions, namely $\Gap(\vec t, \tilde{\vec{\theta}}) \leq \varepsilon$, where $\varepsilon \in \{10^{-5}, 10^{-7}\}$. For the KL-penalized problem, we choose $\lambda \in \{10^{0}, 10^{-1}\}$, and for the $\ell_2$-penalized problem, we select $\lambda \in \{10^{-1}, 10^{-2}\}$.

\section{CONCLUSION}
In this paper, we propose the first Safe Screening method specifically tailored for the UOT, leveraging its unique structure. We introduce a series of novel methods for projection and safe region construction, designed to improve performance without adding significant computational burden. For the KL-penalized UOT, we present a new theoretical analysis framework, from which we derive optimal bounds and prove the property of local blockwise strong convexity. We illustrate the substantial improvements of these methods through numerical experiments. With its capacity to exploit the gradually increasing sparsity, our method holds considerable potential for efficiently solving large-scale problems.

\clearpage
\bibliographystyle{plain}
\bibliography{ref}
\clearpage

\appendix

\section{Notations}
\subsection{The example of the index matrix}

As we illustrated in {\bf Section~\ref{Sec:2.1}}, the index matrices $\mat{N}$ and $\mat{M}$ follow a distinct format. $\mat{N}$ signifies the indices tied to the computation of each row sum for the matrix $\mat{T}$ when presented in vector form, while $\mat{M}$ embodies the column sums. Let us illustrate this with an example where $m=n=3$:\begin{equation*}
\begin{split}
\mat{N}&=\begin{pmatrix}
1&1&1& 0& 0& 0& 0& 0&0\\
0 & 0& 0&1&1&1& 0& 0&0\\
0 & 0& 0& 0& 0& 0&1&1&1\\
\end{pmatrix},\\
\mat{M}&=\begin{pmatrix}
1& 0& 0&1& 0& 0&1& 0&0\\
0&1& 0& 0&1& 0& 0&1&0\\
0& 0&1& 0& 0&1& 0& 0&1\\
\end{pmatrix}.
\end{split}
\end{equation*}

\section{PROOFS}

\subsection{Proofs of {{\bf Theorem~\ref{Thm:Gap_all}} and {\bf Theorem~\ref{Thm:Gap_ellipse}}}}
\begin{proof}
We replace $f$ in \textbf{Definitions}~\ref{Def:SC} with the dual function $\mat D(\theta)$ since
\begin{equation}
\begin{split}
\mat D(\vec y) \leq \mat D(\vec x) +\nabla \mat D(\vec x)^{T}(\vec y-\vec x) - \frac{L}{2}\|\vec y -\vec x\|^{2}
\end{split}
\end{equation}
By choosing $y = \vec{\theta}$ and $x = \hat{\vec{\theta}}$, assuming $\hat{\vec{\theta}}$ is the optimum for the convex function $\mat D$ with $\nabla \mat D(\hat{\vec{\theta}}) = \vec{0}$, we have\begin{equation}
\begin{split}
\mat D({\vec \theta}) &\leq \mat D(\hat{\vec \theta}) - \frac{L}{2}\|\vec \theta -\hat{\vec \theta}\|^{2}\\
\Longleftrightarrow \quad \quad \|\vec \theta -\hat{\vec \theta}\|^{2} & \leq \frac{2}{L}(\mat D(\vec \theta) - \mat D({\hat{\vec \theta}}))\\
\Longleftrightarrow \quad \quad \|\vec \theta -\hat{\vec \theta}\|^{2} & \leq \frac{2}{L}(\mat P(\vec t) - \mat D({\hat{\vec \theta}}))\\
\Longleftrightarrow \quad \quad \|\vec \theta -\hat{\vec \theta}\|^{2} & \leq \frac{2\Gap(\vec t,\vec \theta)}{L}.
\end{split}
\end{equation}
This proves \textbf{Theorem~\ref{Thm:Gap_all}}. As for the proof of \textbf{Theorem~\ref{Thm:Gap_ellipse}}, we can deduce it in the same way by replacing $\frac{L}{2}\| \vec{\theta} - \hat{\vec{\theta}}\|^2_2$ with $\frac{1}{2}(\hat{\vec{\theta}} - \vec{\theta})^T \mathbf{L}(\hat{\vec{\theta}} - \vec{\theta})$ in \textbf{Definitions}~\ref{Def:BSC}.
\end{proof}

\subsection{Proofs of the dual property in {\bf Table~\ref{table:sc}}}

Considering the conjugate function of the UOT function with different penalties:
\begin{proof}
For the $\ell_2$ penalty $h(\vec t) = \|\vec t-\vec y\|_2^2$:
\begin{eqnarray}
h^*(\vec \theta) &=&\max _{\forall p, t_p \geq 0} \ (\vec \theta^{T}\vec t- \|\vec t - \vec y\|^2_2
)\ =\ \frac{1}{2}\|\vec\theta\|_2^2 +\vec y^{T}\vec\theta\notag\\
\end{eqnarray}

For the KL penalty, we have:
\begin{eqnarray}
h^*(\vec \theta) &=&\max _{\forall p, t_p \geq 0} \ (\vec \theta^{T}\vec t-(\vec t +\epsilon)\ln \frac{\vec t + \epsilon}{\vec y}+(\vec t +\epsilon)^{T} \one-\vec y^{T} \one)
\ =\ \vec y^{T} e^{\vec \theta} -\vec y^{T} \one -\epsilon\vec\theta^{T}\one\notag\\
\end{eqnarray}

For the TV penalty, we have:
\begin{eqnarray}
h^*(\vec \theta) &=&\max _{\forall p, t_p \geq 0} \ (\vec \theta^{T}\vec t- \|\vec t - \vec y\|_1)
\ =\ \left\{
\begin{aligned}
      \vec{\theta}^{T}\vec{y},\quad& \forall p,\text{ if }\  |\theta_p| < 1, \\
     +\infty,\quad & \text{otherwise}.
\end{aligned}
\right.\notag\\
\end{eqnarray}
As for $g(x) = \vec c^{T}\vec t$, which does not change with the choice of the penalty function, we have:
\begin{eqnarray}
g^*(\vec\theta) &=&\max _{\forall p, t_p \geq 0} (\vec \theta^{T} \vec t-\lambda \vec c^{T} \vec t)
\ =\ \left\{
\begin{aligned}
0 ,\quad& \forall p \in [n+m], \theta_p \leq \lambda c_p\\
\infty, \quad& \text{otherwise}.
\end{aligned}
\right.
\label{eq_appendix:uotdual}
\end{eqnarray}

The $h^*$ function works as the main dual function and the $g^*$ function works as the constraints. By adding the dual function with the appropriate variable, we can use the Fenchel-Rockafellar Duality to obtain the dual form. During the duality computation, we can derive the primal-dual relationship for $\hat{\vec{\theta}}$ and $\hat{\vec{t}}$, which is used for computing $\vec{\theta}^{k}$ from $\vec{t}^{k}$.\end{proof}

\subsection{Proofs of Theorem.~\ref{Thm:kl_local}}
\begin{proof}
\label{P:bound}
For any $\vec \theta$, we have
\begin{equation*}
\mat P(\hat{\vec t}) = \mat D(\hat {\vec{\theta}}).
\end{equation*}
Substituting $\hat{\vec{\theta}}$ into the dual function yields
\begin{equation*}
\mat P({\vec t}) \geq \mat P(\hat{\vec t}) = \mat D(\hat{\vec\theta}) \geq \mat D(\vec\theta).
\end{equation*}
Now let's consider a single dual variable $\theta_q$:
\begin{equation*}
\mat D(\vec\theta) =\sum_{i\neq q}^{n+m} \mat D^{i}(\theta_i) + \mat D^{q}(\theta_q),
\end{equation*}
where $\mat{D}^{q}(\theta_q) = -y_q e^{-\theta_q} + y_q - \epsilon \theta_q$. We can view $\sum_{i \neq q}^{n+m} \mat{D}^{i}(\theta_i)$ as the dual problem of a new UOT problem $\mat{P}^q(\vec{t})$. In essence, by excluding the element $q$-th element from the dual problem, the corresponding primal problem remains a UOT problem without the element $y_q$. This is akin to modifying the original $\mat{P}$ problem by setting $y_q = 0$, rendering the corresponding $\mat{X}_q$ invalid, along with the indices $t_i, c_i, \text{ that } \forall i \in I_q$. We denote this new $\vec{y}$ as $\grave{\vec{y}}^q$ and establish new $\grave{\vec{c}}^{q}, \grave{\vec{t}}^{q}$, and new $\grave{\mat{X}}^q$ by deleting $y_q$ corresponding elements. Consequently, we derive the corresponding Primal problem for $\sum{i \neq q}^{n+m} \mat{D}^{i}(\theta_i)$, where $\vec{\theta} \in \mathcal{R}^D$:
\begin{equation*}
\mat P^{q}(\grave{\vec t}^{q}) = \lambda {\vec\grave{\vec c}^{qT}}\grave{\vec t}^{q} + \mat D_\phi(\grave{\mat{X}}^q\grave{\vec t}^{q},\grave{\vec y}^q)
\end{equation*}
Then we have
\begin{equation}
\begin{split}
\mat D(\vec\theta) &\leq \mat D(\hat{\vec\theta})\\
& = \sum_{i\neq q}^{n+m} \mat D^{i}(\hat{\theta}_i) + \mat D^{q}(\hat{\theta}_q)\\
& \leq \mat P^{q}({\hat{\grave{\vec t}}^{q}}) + \mat D^{q}(\hat{\theta}_q)\\
& \leq \mat P^{q}(\grave{\vec t}^{q}) + \mat D^{q}(\hat{\theta}_q)
\end{split}
\end{equation}
where $\hat{{\grave{\vec{t}}}}^q$ is the optimum for $\mat P^{q}$. Setting $K^{q}= \mat D(\vec\theta) - \mat P^{q}({\vec{t}})$, we have
\begin{equation*}
\mat D^{q}(\hat\theta_q) > K^q.
\end{equation*}
Then, this implies
\begin{equation*}
- y_q e^{-\hat\theta_q} - \hat\theta_q \epsilon + y_q > K^q.
\end{equation*}
Noting that $\mat D^{q}(\theta_q)$ is an increasing function, we have
\begin{equation*}
(\epsilon - y_p)e^{-\hat\theta_q} - \epsilon > \mat D^{q}(\hat\theta_q)- y_q > K^q - y_q.
\end{equation*}
Finally, we obtain
\begin{equation*}
\hat\theta_p > \ln \left(\frac{\epsilon - y_p}{K_q - y_q +\epsilon}\right) = \operatorname{Low}(\vec \theta, q).
\end{equation*}

We now have the lower bound for $\hat{\theta}_q$. As $\hat\theta_p$ has a lower bound, we can get its upper bound based on the constraints $\vec x_p^{T}\vec\theta < \lambda c_p$, if the $p$-th elements of the flattened vector are in row q or column q-n. Hence, we derive an upper bound for $\hat\theta_q <\min_{p \in I_q} (\lambda c_p - \theta_i)< \min_{p \in I_q} (\lambda c_p -\operatorname{Low}(\vec \theta, i))$from all the relevant constraints. If $0< q \leq n$, then $i \in \{n+1,n+2, ..., n+m \}$, otherwise, if $n< q $, then $i \in \{1,2, ..., n \}$. This completes the proof.
\end{proof}

\subsection{Proof of Theorem.~\ref{Thm:AreaScreeninguOT}}
\begin{proof}
The proposed CTP method constructs a region $\mathcal{R}^{Re}_p$ for every element $t_p$ and makes sure it that includes $\mathcal{R}^{C}\cap \mathcal{R}^{D}$, and is also inside the region $\mathcal{R}^{S}$ proposed in (\cite{Yamada_NIPS_2021}), which is
\begin{equation}
\begin{split}
\mathcal{R}^{S}(\vec{\theta}, \vec{t}):= \left\{\vec{\theta} \ \left|\
\begin{aligned}
&\sum_{p\in I_p}(\vec{x}_{p}^{T}\vec{\theta} - \lambda {c}_p)t_p\leq 0, \sum_{p\in {I}^{C}_p}(\vec{x}_{p} ^{T}\vec{\theta}- \lambda {c}_p)t_p\leq 0, \\
& \vec \theta \in \mathcal{R}^{G}.
\end{aligned}
\right.
\right\}.
\end{split}
\end{equation}
The dual feasible region consists of two parts: one is the circular or elliptical region $\mathcal{R}^{G}$ in {\bf Theorem \ref{Thm:Gap_all}} and {\bf Theorem \ref{Thm:Gap_ellipse}}, which remains unchanged in the construction of the dual feasible area in Sa and CTP. Therefore, it does not require any further consideration, and we focus on the $\mathcal{R}^{D}$ as follows:
\begin{eqnarray*}
\vec \theta \in \mathcal{R}^{D} &\Longrightarrow& \vec x_p^{T}\vec\theta - \lambda c_p \leq 0, \ \forall p \in [nm]\\
&\Longrightarrow& \sum_{p \in I_p \subset [nm]} ( \vec x_p^{T}\vec\theta - \lambda c_p)t_p \leq 0, \ (t_p \geq 0)\\
&\Longrightarrow & \vec \theta \in \mathcal{R}^{Re}_p \\
\end{eqnarray*}
By proving that CTP is a relaxation of $\mathcal{R}^{D}$, since $\vec \theta \in \mathcal{R}^{G}$, we know that the optimum $\hat{\vec \theta} \in \mathcal{R}^{S} = \mathcal{R}^{G}\cap \mathcal{R}^{Re}_p$. This completes the proof.
\end{proof}

\subsection{Proof of Maximizing On $\mathcal{R}^{Re}_{p} \in \mathcal{R}^{G}$}

\begin{proof}
For any $p \in [nm]$, we have:
\begin{equation*}
\begin{split}
\vec{x}_p^T\tilde{\vec{\theta}} &= \tilde{{\alpha}}_{u} + \tilde{{\beta}}_v \\
&= {{\alpha}}_{u} + {{\beta}}_v - \max_{0\leq j< n} \frac{{{\alpha}}_u +{{\beta}}_j - \lambda {c}_{un+j}}{2} - \max_{0 \leq i < m} \frac{{{\alpha}}_i +{{\beta}}_v - \lambda {c}_{in+v}}{2}\\
&\leq{{\alpha}}_{u} + {{\beta}}_v - \frac{{{\alpha}}_u +{{\beta}}_v - \lambda {c}_{p}}{2} - \frac{{{\alpha}}_u +{{\beta}}_v - \lambda {c}_{p}}{2}\\
&\leq \lambda {c}_{p}.
\end{split}
\end{equation*}
This proves that $\forall p$, we have $\vec{x}_p^T\tilde{\vec{\theta}}\leq \lambda {c}_{p} $, then $\tilde{\vec{\theta}} \in \mathcal{R}^{D}$. This completes the proof.
\end{proof}

\subsection{Proof of Maximizing on $\mathcal{R}^{Re}_{p} \cap \mathcal{R}^{G}$ }
\label{P:CTP}
\begin{proof}

By rewriting $\vec \theta := \vec \theta_o + \vec q$, and $\eta,\mu,\nu~(\geq 0)$ are Lagrange multipliers. Optimizing on the intersection of region.~(\ref{Eq:FinalRS}) $\mathcal{R}^{Re}_{p}$ with $\mathcal{R}^{G}$, whose Lagrangian function is given as
\begin{equation*}
\begin{split}
\min_{\vec{q}} \max_{\eta,\mu,\nu \geq 0}
\Biggl\{ L(\vec{q},\eta,\mu,\nu) :=&
- \vec x_p^{T}\vec \theta_o - {q}_{u} - {q}_{u+v+1} + \eta( \vec{q}^T\mat L\vec{q} - 1)\\
&+\mu\left( \left(\sum_{l\in I_p}\vec{x}_{l}\right)^T(\vec{q}+\vec\theta_o) - \lambda\sum_{l\in I_p}c_lt_l\right) + \nu\left( \left(\sum_{l\in I^{C}_p}\vec{x}_{l}\right)^T(\vec{q}+\vec{\theta}_o) -\lambda\sum_{l\in I^{C}_p}c_lt_l\right)
\Biggr\},
\end{split}
\end{equation*}
is the Lagrangian optimization problem of
\begin{equation}
\label{Eq:LagrangianOptProb}
\max_{\vec{\theta} \in \mathcal{R}^{S}_{p}}{ {\theta}_{u} +{\theta}_{u+v+1} },\quad
\forall p \in [nm],
\end{equation}
where $ u = i\mid n,\ v = i\mod n$.

Noticing that, when $\mathcal{R}^{G}$ is a circle, it is just one specific condition by replacing $\mat L$ to a constant L. Considering the center of the circle $\mathcal{R}^{C}$ as $\vec{\theta}_o$, we define $\vec{\theta} := \vec{\theta}_{o} + \vec q$. Also, defining $i_1:=u$, and $i_2:=u+v+1$, since ${\vec x_p^{T} \vec{\theta}_{o}}$ is a constant, Eq.(\ref{Eq:LagrangianOptProb}) can be simplified to
\[
\min_{(\vec{\theta_o} + \vec q )\in \mathcal{R}^{S}_{p}}{- ( {q}_{i_1} +{q}_{i_2} )}.
\]
Here, we define the followings:
\[
\vec{v}:=\sum_{l\in I_p}\vec{x}_{l}, \quad
e_v :=\lambda\sum_{l\in I_p}c_lt_l,\quad
\vec{w}:=\sum_{l\in I^{C}_p}\vec{x}_{l},\quad
e_w :=\lambda\sum_{l\in I^{C}_p}c_lt_l.
\]

Then, we re-write the Lagrangian function as
\begin{equation}
\label{Eq:min_max_2}
\min_{\vec{q}} \max_{\eta,\mu,\nu \geq 0}
\biggl\{ L(\vec{q},\eta,\mu,\nu)
:= - {{q}_{i_1} - {q}_{i_2} + \eta( \vec{q}^T\mat L\vec{q} - 1)+\mu( \vec{v}^T\vec{q} - e_v ) + \nu( \vec{w}^T\vec{q} - e_w )}
\biggr\}.
\end{equation}
Similarly, considering the derivative with respect to ${q}_i$ yields
\begin{equation*}
\begin{split}
\frac{\partial L}{\partial {q}_i} = \left\{
\begin{aligned}
-1 + 2\eta {q}_iL_{i} +\mu v_i + \nu w_i, \quad& i = i_1, i_2,\\
2\eta {q}_iL_{i} +\mu v_i + \nu w_i, \quad& i \neq i_1, i_2.
\end{aligned}
\right.
\end{split}
\label{eq:lang1}
\end{equation*}
Consequently, we obtain
\begin{equation*}
\begin{split}
q_i^{*} = \left\{
\begin{aligned}
\frac{1- \mu v_i - \nu w_i}{2\eta L_{i}}, \quad& i = i_1, i_2,\\
-\frac{\mu v_i + \nu w_i}{2\eta L_{i}}, \quad& i \neq i_1, i_2.
\end{aligned}
\right.
\end{split}
\end{equation*}

Consequently, plugging this into Eq~.(\ref{Eq:min_max_2}), we obtain the following Lagrangian dual problem:
\begin{equation*}
\max_{\eta,\mu,\nu\geq0}
\Biggl\{ L(\eta,\mu,\nu) := \frac{\mu v_{i_1} + \nu w_{i_1}-1}{2\eta L_{i}} +\frac{\mu v_{i_2} + \nu w_{i_2}-1}{2\eta L_{i}}+ \eta(({\vec{q}^{*}})^T\mat L\vec{q}^{*}-1 )+\mu( \vec{v}^T\vec{q}^{*} - e_v ) + \nu( \vec{w}^T\vec{q}^{*} - e_w )
\Biggr\}.
\end{equation*}
From the KKT conditions, we know that if
\begin{eqnarray*}
\eta (({\vec{q}^{*}})^T\mat L\vec{q}^{*} -1) &=& 0,\\
\mu( \vec v^T\vec{q}^{*} - e_v)&=& 0,\\
\nu(\vec w^T\vec{q}^{*} - e_w) &=&0.
\end{eqnarray*}
We set $\eta^{*}, \mu^{*}, \nu^{*}$ as the solution of the equations, which are also the solutions to the dual problem. Firstly, we assume that $\eta^{*}, \mu^{*}, \nu^{*} \neq 0$, then the solution is equal to computing the following equations:
\begin{eqnarray*}
(\frac{1-\mu v_{i_1}-\nu w_{i_1}}{\sqrt{L_{i_1}}})^2 + (\frac{1-\mu v_{i_2}-\nu w_{i_2}}{\sqrt{L_{i_2}}})^2 + \sum^{n+m}_{i\neq i_1,i_2}(\frac{v_i\mu+w_i\nu}{\sqrt{L_i}})^2 - 4\eta^2 &=& 0, \\
\frac{v_{i_1}-\mu v_{i_1}^2-\nu w_{i_1}v_{i_1}}{L_{i_1}} + \frac{v_{i_2}-\mu v_{i_2}^2-\nu w_{i_2}v_{i_2}}{L_{i_2}} - \sum^{n+m}_{i\neq i_1,i_2}(\frac{v_i^2\mu +w_i v_i\nu}{L_{i}}) - 2\eta {e_v} &=& 0, \\
\frac{w_{i_1}-\nu w_{i_1}^2-\mu w_{i_1}v_{i_1}}{L_{i_1}} + \frac{w_{i_2}-\nu w_{i_2}^2-\mu w_{i_2}v_{i_2}}{L_{i_2}} - \sum^{n+m}_{i\neq i_1,i_2}(\frac{w_i^2\nu +w_i v_i\mu}{L_{i}}) - 2\eta {e_w} &=& 0.
\end{eqnarray*}
Set $\vec l = \operatorname{vec}({\mat L})$, and $\tilde{v} = v/l$ then they are rearranged into
\begin{eqnarray*}
(\frac{1}{L_{i_1}}+\frac{1}{L_{i_2}})-2\mu (\tilde{v_{i_1}}+\tilde{v_{i_2}})-2\nu(\tilde{w_{i_1}}+\tilde{w_{i_2}})+ \tilde{\vec{v}}^{T}\vec{v}\mu^2+ \tilde{\vec{w}}^{T}\vec{w}\nu^2+2\mu\nu \tilde{\vec{v}}^T{\vec{w}}- 4\eta^2 &=& 0,\\
(\tilde{v_{i_1}}+ \tilde{v_{i_1}}) - \tilde{\vec{v}}^{T}\vec{v}\mu - \tilde{\vec{v}}^T{\vec{w}}\nu - 2\eta {e_v} &=& 0, \\
(\tilde{w_{i_1}}+ \tilde{w_{i_2}}) - \tilde{\vec{w}}^{T}\vec{w}\nu - \tilde{\vec{v}}^T{\vec{w}} \mu - 2\eta {e_w} &=& 0.
\end{eqnarray*}
From these results, we obtain
\begin{eqnarray*}
\mu &=& \frac{2( e_w\tilde{\vec{v}}^T{\vec{w}} - e_v\tilde{\vec{w}}^{T}\vec{w} )\eta + (\tilde{v}_{i_1}+\tilde{v}_{i_2})\tilde{\vec{w}}^{T}\vec{w} - (\tilde{w}_{i_1} + \tilde{w}_{i_2}) (\tilde{\vec{v}}^T\tilde{\vec{w}})}{ \tilde{\vec{v}}^{T}\vec{v} \tilde{\vec{w}}^{T}\vec{w} -(\tilde{\vec{v}}^T{\vec{w}})^2},\\
\nu &=&\frac{2( e_v\tilde{\vec{v}}^T{\vec{w}} - e_w\tilde{\vec{v}}^{T}\vec{v} )\eta + (\tilde{w}_{i_1}+\tilde{w}_{i_2}) \tilde{\vec{v}}^{T}\vec{v} - (\tilde{v}_{i_1} + \tilde{v}_{i_2}) (\tilde{\vec{v}}^T\tilde{\vec{w}})}{ \tilde{\vec{v}}^{T}\vec{v} \tilde{\vec{w}}^{T}\vec{w} -(\tilde{\vec{v}}^T{\vec{w}})^2}.
\end{eqnarray*}
Denoting 
\begin{eqnarray}
\label{eq:final}
\mu := s_1 \eta + s_2, &\text{and}& \nu := u_1 \eta + u_2,
\end{eqnarray}
we can solve $\eta$ by the following quadratic equation:
\begin{eqnarray*}
0&=&a\eta^2+b\eta+c,\\
a&=&4 - s_1^2\tilde{\vec{v}}^{T}\vec{v} - u_1^2\tilde{\vec{w}}^{T}\vec{w} -2s_1 u_1\tilde{\vec{v}}^T{\vec{w}},\\
b&=&2(\tilde{v}_{i_1} + \tilde{v}_{i_2})s_1 +2(\tilde{w}_{i_1} + \tilde{w}_{i_2})u_1 - 2s_1s_2 \tilde{\vec{v}}^{T}\vec{v} - 2u_1u_2\tilde{\vec{w}}^{T}\vec{w} - 2(s_1u_2+s_2u_1)\tilde{\vec{v}}^T{\vec{w}},\\
c&=&2(\tilde{v}_{i_1} + \tilde{v}_{i_2})s_2 +2(\tilde{w}_{i_1} + \tilde{w}_{i_2})u_2 -s_2^2\tilde{\vec{v}}^{T}\vec{v} -u_2^2\tilde{\vec{w}}^{T}\vec{w} - 2s_2u_2\tilde{\vec{v}}^T{\vec{w}} -(\frac{1}{L_{i_1}}+\frac{1}{L_{i_2}}).
\end{eqnarray*}

Putting these equations back into Eq.(\ref{eq:final}), we obtain $\mu$ and $\nu$.

If the solution satisfies the constraints $\eta^{*}, \mu^{*}, \nu^{*} > 0$, then it is exactly the solution for the optimization problem.
However, if one of the dual variables is less than 0, the problem would degenerate into a simpler question as one of the constraints is not activated.

1) If only the circle is activated:, then $\min_{\vec{\theta} \in \mathcal{R}^{S}_{I}}{- ( {q}_{i_1} +{q}_{i_2} )} = \frac{-2}{\sqrt{L_1 + L_2}}$.

2) If only one plane and one circle are activated, we can optimize on a spherical cap or dome, the solution is the same with (\cite{Yamada_NIPS_2021}). Or just setting the $\nu = 0$ in our formula. As the hyperplane $\vec{w}^T\vec{q} - e_w $ is parallel to the optimization direction, it is impossible to be the activated plane with the circle.
Thanks to the distinctive structure of the index matrix of the UOT problem, the summation involved in solving this Lagrangian problem, with a complexity of $\mathcal{O}(m)$ or $\mathcal{O}(n)$, can be precomputed. Subsequently, this precomputed summation can be reused for elements within the same row or column, thereby maintaining the overall complexity of the screening method at $O(nm)$. This complexity is consistent with other screening methods for the UOT problem.
\end{proof}

\subsection{Computational complexity}
The complexity of Safe Screening for $\ell_2$-penalized UOT primarily consists of three components: 1) duality and projection, 2) computation of the maximum value within the safe region. 

Regarding the duality and projection part, we have discussed it in the main text, where the complexity of the Shifting Projection is nearly identical to that of Residuals Rescaling. As for the computation of the maximum value, we have addressed it in \textbf{Section~\ref{P:CTP}}. Although the CTP method involves additional computations compared to the Gap method, it shares the same complexity as the Dome method and the projection process.

However, for KL-penalized OT, an additional step is introduced in the Safe Screening procedure: constructing a $\mathcal{R}^{B}$ with local strong concavity, which necessitates the computation of bounds for optimizing the objective function. Based on the insights from the proof presented in \textbf{Section~\ref{P:bound}}, we observe that the computation of bounds is similar to finding the optimum in the CTP method. Due to the specific structure of the UOT problem, variables related to row and column summations can be reused, resulting in a complexity of $O(nm)$. Thus, there is no increase in computational complexity in terms of the order of magnitude.

\subsection{Proof for Optimal Stepsize of FISTA in Section~\ref{Sec:EX1}}
The Fast Iterative Shrinkage-Thresholding Algorithm (FISTA)\cite{Beck_2009_SIAMIS} is a popular solver for regularized problems and can be applied to $\ell_1$ regularized problems. It is a variant of the Mirror Descent algorithm, specifically when the smoothing part in the proximal operator is chosen as $\|\vec{x}\|_2^2$. The theoretical optimal step size for the Mirror Descent method is $1/L$, where $L$ is the constant of $L$-relatively smoothness for the function $H(\vec{t}) = h(\mat{X} \vec{t})$ \cite{doi:10.1137/16M1099546}. 

In the context of applying FISTA to the UOT problem, where the proximal function is $\phi(\vec{t}) = \|\vec{t}\|_2^2$, we aim to compute the minimum constant $L$ that ensures the relative smoothness of the UOT problem. It is important to note that if a function $H$ is $L$-relatively smooth with respect to function $\phi$, then the function $L\phi - H$ is convex, or equivalently, $D_H(\vec{x}, \vec{y}) < LD_{\phi}(\vec{x}, \vec{y})$ \cite{bauschke2017descent}

\begin{proof}
We prove that $L\phi - H$ is convex by proving its positive semi-definiteness, i.e., $\vec{d}^T \nabla^2(L\phi - H)\vec{d} \succeq 0$ for any $\vec{d}\in{\mathbb{R}^{nm}}$. Here, we denote that the $i$-th row vector of $\mat{X}$, i.e., $\vec{x}^i = [X_{i,1},X_{i,2},X_{i,1},...,X_{i,nm}]$. Then, because $\nabla^2 H(\vec{t}) = \mat X^{T}\mat X$ where $\mat{X} \in \R^{(n+m) \times nm}$, we have
\begin{eqnarray*}
\vec{d}^{T}\nabla^2 H(\vec t)\vec{d} &=& \sum_{i=1}^{n+m}{(\vec{d}^{T}\vec{x}^i)^2}\\
&=& 2\sum_{j=1}^{nm} d_j^{2} + 2\sum_{p\in I_p}^{mn}{d_u d_v}\\
& \leq& 2\sum_{j=1}^{nm} d_j^{2} + (n+m - 2)\sum_{j}^{nm}{d_j^2}\\
&=& (n+m) \sum_{j=1}^{nm}d_j^{2}\\
&\leq& L \sum_{j}^{nm}{d_{j}^{2}}\\
&=& L \vec{d}^{T} \nabla^2 \phi(\vec t)\vec{d}.
\end{eqnarray*}
Thus, we have proven the convexity of the $L\phi - H$ as $L\geq n+m$, and the theoretical stepsize is $\frac{1}{n+m}$.
As for the MM algorithm in {\cite{Chapel_NeurIPS_2021}} following the same method, we can get a fixed stepsize as $\frac{1}{2}$, which is the same as the parameter in the paper.
\end{proof}

\section{ADDITIONAL EXPERIMENTS}

This section provides some additional experimental results.

\subsection {Comparison on coordinate descent (CD) method} 
We have organized our experiments on CD with $\ell_2$-penalized UOT, which has the same parameters setting with {\bf Para~\ref{para:ex3}}. Our experiments show great improvement, especially on the MNIST dataset. The outcome is illustrated in {\bf Table {\ref{table:CD}}}.

\begin{table}[h]
\caption{Speed-up ratios of the Safe Screening: CD for the $\ell_2$ penalty.}
\label{table:CD}
\centering
\begin{tabular}{l|c|cc|cc}
\toprule
Dataset & Method & \multicolumn{2}{c}{$\lambda=10^{-1}$} & \multicolumn{2}{c}{$\lambda=10^{-2}$} \\
\cmidrule{3-6}
& & $10^{-5}$ & $10^{-7}$ & $10^{-5}$ & $10^{-7}$ \\
\midrule
MNIST & Gap & 1.55 & 3.67 & 0.92 & 1.26 \\
& Sa & 1.78 & 4.34 & 0.94 & 1.42 \\
& Sa-CTP & $\bm{2.17}$ & $\bm{5.12}$ & $\bm{1.03}$ & $\bm{2.09}$ \\
\midrule
Gauss & Gap & $\bm {1.52}$ & $\bm {2.83}$ & $\bm {1.03}$ & 1.77 \\
& Sa & 1.46 & 2.56 & 1.00 & 2.26 \\
& Sa-CTP & 1.50 & 2.33 & 0.97 & $\bm{2.87}$ \\
\bottomrule
\end{tabular}
\end{table}
\subsection{Comparison on BFGS method}
We have organized our experiments on L-BFGS-B solver \cite{DBLP:journals/toms/MoralesN11} with both $\ell_2$-penalized and KL-penalized UOT, which has the same parameters setting with {\bf Para~\ref{para:ex3}}. The outcome is illustrated in {\bf Table {\ref{table:BFGS}}}. 

It is worth noting that the original code for the L-BFGS-B algorithm is encapsulated in Fortran, making direct modifications challenging for us. Consequently, we choose to call the L-BFGS-B Fortran code via {\it mex} and restart the algorithm after the Safe Screening process. This approach, albeit interrupting the update of the Hessian matrix, could potentially impact the convergence speed of the algorithm. In our comparison, both methods, with and without screening, utilized this interrupt-and-restart optimization process. This was undertaken to ensure the consistency of our experimentation. We only add BFGS experiments on KL-penalized UOT with $\lambda = 1$ because the convergence is too slow.

\begin{table}[h]
\caption{Speed-up ratios of the Safe Screening: L-BFGS-B for $\ell_2$ penalty and KL penalty.} 
\label{table:BFGS}
\centering
\begin{tabular}{l|c|cc|cc|c|cc}
\toprule
Dataset & Method & \multicolumn{4}{c|}{$\ell_2$} & Method & \multicolumn{2}{c}{KL} \\
& & \multicolumn{2}{c}{$\lambda=10^{-1}$} & \multicolumn{2}{c|}{$\lambda=10^{-2}$} & & \multicolumn{2}{c}{$\lambda=1$} \\
\cmidrule{3-6} \cmidrule{8-9}
& & $10^{-5}$ & $10^{-7}$ & $10^{-5}$ & $10^{-7}$ & & $ 10^{-5}$ & $10^{-7}$ \\
\midrule
MNIST & Gap & 1.65 & 3.00 & 1.01 & 1.24 & Gap & 1.01 & 1.04 \\
& Sa & 1.92 & 3.47 & 1.01 & 1.37 & Gap-CTP& 0.99 & 1.07 \\
& \!Sa-CTP\! & $\bm{2.25}$ & $\bm{4.18}$ & $\bm{1.11}$ & $\bm{1.66}$ & Ell & 1.08 & 1.23 \\
& & & & & & \!\!Ell-CTP\!\! & $\bm{1.09}$ & ${\bm{1.24}}$ \\
\midrule
Gauss & Gap & $\bm {1.39}$ & 2.25 & 1.03 & 1.88 & Gap & 1.04 & 1.01 \\
& Sa & 1.32 & ${\bm {2.87}}$ & ${\bm {1.03}}$ & 2.18 & Gap-CPT & 1.04 & 1.01 \\
& Sa-CTP & 1.19 & 2.65 & 0.97 & $\bm{2.64}$ & Ell & 1.75 & 1.87 \\
& & & & & & \!\!Ell-CTP\!\! & $\bm{1.81}$ & ${\bm {2.18}}$ \\
\bottomrule
\end{tabular}
\end{table}

\section{Discussion and Perspectives}

\subsection{Large-scale OT Problem and Sparsity}
The sparsity of the OT and UOT problems is $1 - 1/(n+m)$ \cite{doi:10.1137/16M1106018}. This indicates that a sparse-aware algorithm could potentially achieve a speed-up factor of $n+m$. Employing Safe Screening for OT and UOT problems is therefore promising for large-scale problems.

\subsection{Aggressive Screening}
Our Shifting Projection and CTP methods exploit the index matrix structure characteristics, and the Ellipse method utilizes the anisotropy of the dual variables in the dual space. However, there is potential for the development of more aggressive Safe Screening methods. For instance, the locally strongly concave function we proved on $\mathcal{R}^B$ indicates that the Box bound is quite loose compared to the real dual optimum. This suggests that an improved method could enhance the screening performance on KL-penalized UOT by several orders of magnitude. Additionally, in real applications, we do not require a completely \textit{safe} screening, which means that we could utilize more relaxed bounds to screen more elements, thereby achieving a faster algorithm with a minor sacrifice in accuracy.

\subsection{Limitations and Shortcomings}
Due to the specific structure of UOT, we can calculate its theoretical Lipschitz constant to compute an theoretically optimal stepsize. However, Safe Screening changes the dimension and structure of the problem, causing the optimal stepsize to change dynamically. Moreover, an effective Safe Screening method would identify theoretically zero elements before they reach zero. This behavior could detrimentally impact the error (like dual gap), and temporarily slow down the screening process. Such oscillation is difficult to manage. We prevent the algorithm from screening elements larger than zero to maintain the stability of our algorithm, facilitating a straightforward comparison of the speed-up ratio. This limitation also inhibits the potential of Safe Screening. Further research is necessary to facilitate its application.
\end{document}